\documentclass[12pt,oneside]{amsart}
\usepackage{amsxtra}
\usepackage{amsopn}
\usepackage{amsmath,amsthm,amssymb}
\title{ Almost K\"ahler 4-dimensional Lie groups with $J$-invariant Ricci tensor}
\author{Anna Fino}
\thanks{Research partially supported by MURST, GNSAGA (Indam) of Italy and EDGE Research Training Network HPRN-CT-2000-00101.}
\address{Dipartimento di Matematica, Universit\'a di Torino\\
Via Carlo Alberto 10, 10123 Torino, Italy}
\email{annamaria.fino@unito.it}

\newtheorem{teo}{Theorem}[section]

\newtheorem{corol}{Corollary}[section]

\newtheorem{remark}{Remark}[section]

\newcommand{\beq}{\begin{equation}}
\newcommand{\eeq}{\end{equation}}
\newcommand{\bqn}{\begin{eqnarray}}
\newcommand{\eqn}{\end{eqnarray}}
\newcommand{\bqne}{\begin{eqnarray*}}
\newcommand{\eqne}{\end{eqnarray*}}

\newcommand{\R}{{\mathbb R}}

\begin{document}

\maketitle

\begin{abstract} The aim of this paper is to determine left-invariant strictly almost K\"ahler structures on 4-dimensional Lie
groups $(g, J, \Omega)$ such that the Ricci tensor is $J$-invariant.

\end{abstract}
\section{Introduction}

An almost K\" ahler structure on a manifold $M$ of real dimension $n = 2m$ consists of a
symplectic form $\Omega$, an almost complex structure $J$ and a Riemannian metric $g$, 
satisfying the compatibility condition
\begin{equation} \label{compatibility}
\Omega (X, Y) = g (J X, Y),
\end{equation}
for any vector fields $X,Y$  on $M$. 
If $J$ is integrable, then the triple $(g, J, \Omega)$ is a K\" ahler structure on $M$. In real dimension 2, the notions of
almost K\" ahler and K\" ahler structure coincide, but this does not hold in higher dimensions. Through the paper,  ``strictly
almost  K\" ahler" will mean that the corresponding almost complex structure is non integrable (or, equivalently, that the
almost K\" ahler structure is non K\" ahler).

Given a symplectic manifold $(M, \Omega)$, there are many $\Omega$-compatible almost K\" ahler structures on $M$, i.e. many 
pairs
$(g, J)$ which satisfy the relation
\eqref{compatibility}. The space ${\mathcal AK} (M, \Omega)$ of all $\Omega$-compatible metrics (or almost complex structures)
is an infinite dimensional, contractible Fr\'echet space. The problem of finding distinguished Riemannian metrics in 
${\mathcal AK} (M, \Omega)$ has been intensively studied (for a survey, see \cite{AD}). In
particular, such a kind of  metric is given by an Einstein one. A still open conjecture by Goldberg affirms that there
exist no Einstein, strictly almost K\"ahler metrics on a compact symplectic manifold \cite{Go}. This conjecture is true if the
scalar curvature is non-negative \cite{Se} and there many positive results in dimension four under different additional
assumptions on the curvature \cite{AA,A1,A2,OS}. Moreover, compactness is important in the above  conjecture since  in
\cite{NP} an example of a
non-compact Ricci flat strictly almost K\"ahler manifold  is given. 

On a compact manifold $M$, the Einstein condition agrees with the Euler-Lagrange equation of the Hilbert functional (integral
of the scalar curvature) acting on the space of all Riemannian metrics on $M$ of a given volume \cite{Be}. 
Blair and Ianus \cite{BI} restricted such  a functional to the space ${\mathcal AK} (M, \Omega)$ and
proved that its critical points are the almost K\"ahler metrics $(g, J)$ whose Ricci tensor $\rho$
is
$J$-invariant, i.e. satisfies
\begin{equation} \label{Ricci}
\rho(J X, J Y) = \rho (X,Y),
\end{equation}
for any vector fields $X,Y$. Therefore, the previous condition is  weaker than the Einstein and K\"ahler condition. Almost
K\"ahler metrics $(g, J)$ satisfying \eqref{Ricci} have been recently studied in \cite{LW} and have been
called \lq \lq harmonic almost K\" ahler structures". 

In this paper we study such harmonic strictly almost K\"ahler structures on simply-connected  4-dimensional Riemannian
homogeneous spaces, in particular on  4-dimensio\-
nal real Lie
groups. In general,  any simply connected  4-dimensional Riemannian homogeneous space is either symmetric or
isometric to a Lie group with a left-invariant metric \cite{BB}. By \cite{MOS} there exists no a compact strictly almost
K\" ahler locally symmetric space.

 The problem of existence of special left-invariant metrics, like
Einstein and anti-self-dual metrics, have been studied respectively  in \cite{Je} and
\cite{DS}, where a classification is obtained.

The almost K\" ahler  structures $(g, J, \Omega)$ that we consider are invariant, in the sense that
$g, J, \Omega$  are left-invariant tensors.  We prove that a simply-connected 4-dimensional  real Lie group with an
invariant harmonic strictly almost K\" ahler structure is necessarily solvable and we  give a  general description of  its Lie
algebra. Moreover, we prove that all these Lie groups are isometric (up to homothety) to the (unique) 4-dimensional proper
3-symmetric space (described by Kowalski \cite{Ko})
 and satisfy the condition \cite{Gr} 
$$
(G2) \quad R_{XYZW} = R_{JX JY Z W} + R_{JX Y JZ W} + R_{JX Y Z JW}.
$$
This space is given by $\R^4$ endowed with
the Riemannian metric
$$
\begin{array} {ll}
g = &(-x + \sqrt{ x^2 + y^2 + 1}) d u^2 +  (x + \sqrt{ x^2 + y^2 + 1}) dv^2 - 2 y du dv\\
{}& +
\lambda^2 (1 + x^2 + y^2)^{-1} [ (1 + y^2) dx^2 + (1 + x^2) dy - 2 xy dx dy)  ], \quad
\lambda > 0.
\end{array}
$$
The underlying homogeneous space $M$ is  $ \R^2 \ltimes SL(2, \R)/SO(2)$
$$
\left \{ \left( \begin{array} {ccc} a&b&u \\
c&d&v\\
0&0&1
\end{array} \right) \backslash \left( \begin{array} {ccc} \cos {t}& - \sin{t}&0\\
\sin{t}&\cos {t}&0\\
0&0&1
\end{array} \right), \quad
\mbox{det} \left( \begin{array} {cc} a&b\\
c&d \end{array} \right) = 1 \right\}
$$
and admits a simply transitive action of a 4-dimensional (real) Lie group $\R^2 \ltimes {\mbox Sol}_2$ \cite{AD}, where
${\mbox Sol}_2$ consists of the upper triangular matrices with positive diagonal entries, \cite{AD,Wa}. Choosing a basis of
the Lie algebra of $\R^2 \ltimes {\mbox Sol}_2$ consisting of   $\{ e_1,
e_2 \}$ in $\R^2$ and 
$$
e_3 = \left( \begin{array} {cc} 1 & 0\\ 0 & -1 \end{array} \right), \quad e_4 = \left( \begin{array} {cc} 0 & 1\\ 0 & 0
\end{array} \right)
$$
in $ {\mathfrak sol}_2 \subset {\mathfrak {sl}}_2 (\R)$, the following are the non zero Lie brackets
$$
[e_1, e_3] = - e_1 = [e_2, e_4], \quad [e_2, e_3] = e_2, \quad [e_3, e_4] = 2 e_4.
$$

\section{Structure equations}

Let $G$ be a real 4-dimensional Lie group with a left-invariant Riemannian metric $g$ (which corresponds to an inner product
$g$ on the Lie algebra ${\mathfrak g}$ of $G$). 
 The curvature tensor $R$ of a left-invariant Riemannian metric is completely
determined by its value at the identity element of $G$. Hence $R$ is an element of the vector space
$$
{\mathcal R} = S^2 (\Lambda^2 {\mathfrak g}^*) \oplus \Lambda^4 {\mathfrak g}^*,
$$
where ${\mathfrak g}$ is the Lie algebra of $G$,
and the Ricci tensor $\rho$ can be identified as the component of $R$ on $S^2 ({\mathfrak g}^*)$.

Then, the curvature tensor of $g$ and its Ricci tensor are completely determined by the structure constants of the Lie algebra
${\mathfrak g}$ and its inner product. Therefore the problem is purely algebraic.

\begin{teo} If $G$ is a simply-connected real 4-dimensional Lie group with an invariant strictly
almost K\"ahler structure $(g, J, \Omega)$ such that the Ricci tensor is $J$-invariant then there exist an orthonormal basis
$\{ e^1, e^2, e^3, e^4 \}$ of left-invariant 1-forms such that
$$
\begin{array} {l}
\Omega  = e^1 \wedge e^2 + e^3 \wedge e^4,\\
J e^1 = e^2, J e^3 = e^4,\\
\end{array}
$$
$$
\begin{array} {l}
d e^1 = - s (e^1 \wedge e^3 - e^2 \wedge e^4) + \frac{(-s^2 + t^2)}{2 t} e^1 \wedge e^4 + t
e^2 \wedge e^3,\\ d e^2 = - \frac{s^2}{t} e^1 \wedge e^3 + s (e^1 \wedge e^4 + e^2 \wedge
e^3) +\frac{(s^2 - t^2)}{2 t} e^2 \wedge e^4,\\
d e^3 = \frac{(t^2 + s^2)} {t} e^3 \wedge e^4, d e^4 = 0,
\end{array}
$$
where $t \neq 0$ and $s$ are real numbers. Relative to the above orthonormal basis, the Ricci tensor has the diagonal form
$$
\left( \begin{array} {cccc}
0&0&0&0\\
0&0&0&0\\  0&0&\frac{ -3(t^2 + s^2)^2}{2 t^2}&0\\
0&0&0&\frac{ -3(t^2 + s^2)^2}{2 t^2} \end{array} \right).
$$
\end{teo}
 
\begin{proof} Since $(G, \Omega)$ is a 4-dimensional Lie group with a left-invariant symplectic
structure $\Omega$, $G$ has to be solvable \cite{LM}. Then the dimension of the commutator ${\mathfrak g}^1 =[
{\mathfrak g}, {\mathfrak g}]$ is less than equal to 3.
One can suppose that there exists 
an orthonormal basis of left-invariant 1-forms
$\{ e^1, e^2, e^3, e^4 \}$ on $G$
such that
$$
\begin{array} {l}
\Omega  = e^1 \wedge e^2 + e^3 \wedge e^4,\\
J e^1 = e^2, J e^3 = e^4,\\
d e^1 = a_1 e^1 \wedge e^2 + a_2 e^1 \wedge e^3 + a_3 e^1 \wedge e^4 +a_4 e^2 \wedge e^3 +a_5 e^2 \wedge e^4 +a_6 e^3 \wedge
e^4,\\
d e^2 = b_1 e^1 \wedge e^2 + b_2 e^1 \wedge e^3 + b_3 e^1 \wedge e^4 +b_4 e^2 \wedge e^3 +b_5 e^2 \wedge e^4 +b_6 e^3 \wedge
e^4,\\
d e^3 = c_1 e^1 \wedge e^2 + c_2 e^1 \wedge e^3 + c_3 e^1 \wedge e^4 +c_4 e^2 \wedge e^3 +c_5 e^2 \wedge e^4 +c_6 e^3 \wedge
e^4,\\
d e^4 = 0,
\end{array}
$$
where $a_i, b_i, c_i$, $i = 1, \ldots, 6$, are arbitrary real numbers. Thus, to describe the 4-dimensional Lie groups
with a strictly almost K\"ahler structure such that the Ricci tensor is $J$-invariant is equivalent to
determine the real numbers $a_i, b_i, c_i$, $i = 1, \ldots, 6$, such that 
\begin{equation} \label{conditions}
\left\{ \begin{array} {l}
d^2 e^i = 0, \quad i = 1, 2, 3, \quad 
d \Omega = 0,\\
\rho(e_1, e_1) = \rho(e_2, e_2), \quad 
\rho(e_1, e_2) = 0,\\
\rho(e_1, e_3) = \rho(e_2, e_4), \quad
\rho(e_1, e_4) = -\rho(e_2, e_3),\\
\rho(e_3, e_3) = \rho(e_4, e_4), \quad
\rho(e_3, e_4)  = 0,
\end{array} \right.
\end{equation}
(where $\{ e_1, \ldots, e_4 \}$ is the dual basis of $\{ e^1, \ldots, e^4 \}$ ) and $J$ is non integrable.
 
The condition $d \Omega =0$ is equivalent to 
\begin{equation} \label{eq1} \left\{
\begin{array} {l}
c_1 = a_3+b_5,\\
  a_2+b_4=0,\\
  b_6-c_2=0,\\
 a_6+c_4=0
\end{array} \right.
\end{equation}
and the conditions $d^2 e^i = 0$, $i = 1, 2, 3$, are:
\begin{equation} \label{eq2} \left\{
\begin{array} {l}
-a_1 b_6-c_6 a_2-b_3 a_4+a_5 b_2+b_6 a_6 = 0,\\
-a_1 a_2-b_1 a_4-b_6 a_4-a_6 a_2=0,\\
a_1 b_5 +a_2 c_5- a_3 a_6-c_3 a_4-b_1 a_5-b_5 a_6=0,\\
a_1 a_6-2 a_2 a_5+a_3 a_4-b_5 a_4-c_6 a_4-a_6^2 = 0,\\
-b_1 b_6-a_3 b_2-c_6 b_2+2 a_2 b_3+b_5 b_2+b_6^2 =0,\\
b_1 a_2-a_1 b_2-a_6 b_2+b_6 a_2=0,\\
b_1 a_3+b_2 c_5-a_1 b_3+a_2 c_3- b_5 b_6-a_3 b_6=0,\\
b_1 a_6-a_5 b_2+b_3 a_4+c_6 a_2-b_6 a_6=0,\\
-2 a_3 b_6-b_5 b_6+a_6 b_3+c_3 a_2+b_2 c_5=0,\\
2 a_3 b_5+b_5^2+a_3^2+ b_6 c_5+ c_3 a_6-a_1 c_3-c_6 a_3-c_6 b_5-b_1 c_5=0,\\
a_3 a_6+2 b_5 a_6-b_6 a_5+c_3 a_4-a_2 c_5=0,\\
a_1 b_6 - a_6 b_1 =0.
\end{array} \right.
\end{equation}
The Ricci tensor $\rho$ of the metric $g$   is $J$-invariant
if and only if
\begin{equation} \label{eq3} \left\{
\begin{array} {l} 
\rho(e_1,e_1) - \rho(e_2,e_2) =
 \frac{1}{2} c_5^2
-c_6 a_3
-b_5 a_4+a_4^2-b_2^2+a_5^2-b_3^2+c_6 b_5\\
 - \frac{3}{2} b_6^2+ \frac{3}{2} a_6^2
+b_5^2-a_3 b_2
-a_3 a_4 -a_3^2-b_5 b_2 - \frac{1}{2} c_3^2-a_1 a_6+b_1 b_6=0,\\
{}\\
\rho(e_1,e_2) =
-b_1 a_6-a_2 a_4
+\frac{3}{2} b_6 a_6
- \frac{1}{2} c_3 c_5 -b_3 b_5+a_2 b_5\\
a_2 a_3 +a_2 b_2 -a_5 a_3 - \frac{1}{2} c_6 a_5- \frac{1}{2} c_6
b_3=0,\\
{}\\
 \rho(e_1,e_3) - \rho(e_2,e_4) =
 a_6 c_6+\frac{3}{2} b_6 a_2- \frac{1} {2} b_6 a_5\\
- \frac{1}{2} a_4 a_6-c_6 c_3+ \frac{1}{2} a_5 c_5+ \frac{1}{2} a_3
a_6- a_6b_5 \\
 +
\frac{1} {2} a_2 c_5-\frac{3} {2} a_6 a_3-
 a_6 b_2\\ 
\frac{1}{2} c_3 a_3- \frac{1}{2} c_3 a_4+a_1 a_3+a_1 a_4+b_1 b_3-b_1 a_2 =0,\\
{}\\
\rho(e_1,e_4)+ \rho(e_2,e_3) =
-\frac{1} {2} a_6 b_3+ \frac{3} {2} a_6 a_2\\
+   a_4 b_6+  \frac{1}
{2} c_5 b_5+  \frac{1} {2} c_5 b_2+  \frac{1} {2} c_3 b_3+  \frac{1} {2} c_3 a_2
-   a_3 b_6  -  \frac{3}
{2} b_6 b_5+ b_6 c_6\\
-a_1 a_2+a_1 a_5-b_1 b_2+b_1 b_5-c_6 c_5+  \frac{1} {2} b_2 b_6 =0,\\
{}\\
\rho(e_3,e_3) - \rho(e_4,e_4) = -
2 a_2^2+c_5^2
-c_6 a_3- \frac{1} {2} a_4^2-a_4 b_2\\- \frac{1}
{2} b_2^2 +a_5 b_3+ \frac{1} {2}a_5^2+ \frac{1} {2} b_3^2-c_6 b_5+b_6 c_5+c_3 a_6-  a_6^2\\
-  b_6^2+\frac{3}{2} b_5^2+a_3 b_5+ \frac{3}{2} a_3^2- b_1 b_6- a_1 a_6+c_3^2 =0,\\
{}\\
\rho(e_3,e_4) =
a_6 c_5- \frac{1} {2} b_1 c_3+ \frac{1} {2}
b_1 a_6-  a_2 a_3+a_6 b_6-b_6 c_3\\
-b_6 a_6+\frac{1} {2} a_1 c_5- \frac{1} {2} a_1 b_6- \frac{1}
{2} b_3 a_4- \frac{1} {2} b_3 b_2\\  a_2 b_5- \frac{1} {2} a_5 a_4- \frac{1} {2} a_5 b_2 =0.
\end{array} \right.
\end{equation}

Moreover, $J$ is integrable if and only 
\begin{equation} \label{eqintegrability} \left\{
\begin{array} {l}
c_4 = - c_3, c_5 = c_2,\\
b_4 = - b_3 + a_2 - a_5,\\
b_5 = b_2 + a_3 + a_4.
\end{array} \right.
\end{equation}

Then, by imposing \eqref{eq1},\eqref{eq2},\eqref{eq3},  we get  the following cases for the relations between the non zero
parameters:
\newline
\smallskip
\noindent
{\sl case 1)}  \hskip 10pt   $b_3 = - a_5$, $c_6$ any real number; \newline
{\sl case 2)} \hskip 10pt $c_6 = 2 b_5, b_3 = - a_5, a_3 = b_5$; \newline
{\sl case 3)} \hskip 10pt $c_4 = - a_6$, $a_1 = \frac{(a_4^2 + a_4 c_6 + a_6^2)}{a_6}, a_3 = - a_4,c_3 = a_6$; \newline
{\sl case 4)} \hskip 10pt $c_6 = -a_4, a_3 = - a_4$, $a_1$ any real value; \newline
{\sl case 5)} \hskip 10pt $c_6 = 2 b_5, a_3 = - b_5, b_2 = - 2 b_5$; \newline
{\sl case 6)} \hskip 10pt $c_6 = 2 b_5, a_3 = - b_5, b_2 =  2 b_5$;\newline
{\sl case 7)} \hskip 10pt $a_4 = - b_2, a_5 = - b_3$;\newline
{\sl case 8)} \hskip 10pt $a_2 = - b_4, a_3 = \frac{b_4^2} {b_5},  a_4 = - \frac{b_4^2} {b_5}, a_5 = - b_4, b_2 = b_5, b_3
= - b_4, c_6 = \frac{(b_4^2 + b_5^2)} {b_5}$;\newline
{\sl case 9)} \hskip 10pt $a_1, b_1, c_6$ any real value; \newline
{\sl case 10)} \hskip 10pt $a_3 = - \frac{1}{2} b_2, b_5 = \frac{1}{2} b_2, c_6 = b_2$;\newline
{\sl case 11)} \hskip 10pt $a_3 =  \frac{1}{2} b_2, b_5 = -\frac{1}{2} b_2, c_6 = -b_2$;
\newline
{\sl case 12)} \hskip 10pt $b_1 = \frac{(- c_6 b_2 + b_2^2 + b_6^2)}{b_6}, b_5 = b_2, c_2 = b_6, c_5 = b_6$; \newline
{\sl case 13)}  \hskip 10pt $a_1 = a_6, b_1 = b_6, c_2 = b_6, c_3 = a_6, c_5 = b_6$, $c_6$ any real value; \newline
{\sl case 14)} \hskip 10pt $a_3 = \frac{1}{2} c_6, a_4 = c_6, b_5 = - \frac{1}{2} c_6$;
\newline
{\sl case 15)} \hskip 10pt $a_3 = \frac{1}{2} c_6, a_4 = -c_6, b_5 = - \frac{1}{2} c_6$;
\newline
{\sl case 16)} \hskip 10pt $a_2 = - a_5, a_3 = \frac{(- a_5^2 + a_4^2)}{2 a_4}, b_2 = - \frac{a_5^2}{a_4}, b_3 = a_5, b_4 =
a_5, b_5 =
\frac{(a_5^2 - a_4^2)}{2 a_4}$,\newline
$c_6 = 
\frac{(a_4^2 + a_5^2)} {a_4}$; \newline
{\sl case 17)}  \hskip 10pt $a_2 =  a_5, a_3 = \frac{ (a_5^2 - a_4^2)}{2 a_4}, b_2 = - \frac{a_5^2}{a_4}, b_3 = a_5, b_4 =
-a_5, b_5 = \frac{(-a_5^2 + a_4^2)}{2 a_4}, \newline
c_6 = 
\frac{(-a_4^2 - a_5^2)} {a_4}$.

Using \eqref{eqintegrability} we can exclude the cases 1), 2), 3), 4), 6), 7), 8), 9), 10), 12),
13), 15), since in all these cases
$J$ is integrable. In the remaining cases the corresponding Lie group has the following structure
equations:
$$
5) \left\{ \begin{array} {l}
d e^1 = - b_5 e^1 \wedge e^4, d e^2 = b_5 (e^2 \wedge e^4 - 2 e^1 \wedge
e^3),\\
d e^3 = 2 b_5 e^3 \wedge e^4, d e^4 = 0, 
\end{array} \right.
$$
(with $\rho(e_1, e_1) = 0$ and $\rho(e_3, e_3) = - 6 b_5^2$)
$$
11) \left\{ \begin{array} {l}
d e^1 =  \frac{1}{2} b_2 e^1 \wedge e^4, d e^2 = b_2 e^1 \wedge e^3 - \frac{1}{2} b_2 e^2 \wedge e^4,\\
d e^3 = -b_2 e^3 \wedge e^4, d e^4 = 0,
\end{array} \right.
$$
(with $\rho(e_1,e_1)=0$ and $\rho(e_3,e_3) = - \frac{3}{2} b_2^2$)
$$
14) \left\{ \begin{array} {l}
d e^1 = \frac{1}{2} c_6 e^1 \wedge e^4 + c_6 e^2 \wedge e^3, d e^2 = - \frac{1}{2} c_6 e^2 \wedge e^4,\\
d e^3 = c_6 e^3 \wedge e^4, d e^4 = 0,
\end{array} \right.
$$
(with $\rho(e_1, e_1) = 0$ and $\rho(e_3,e_3)= - \frac{3} {2} c_6^2$)
$$
16) \left\{
\begin{array} {l}
d e^1 = - a_5 (e^1 \wedge e^3 - e^2 \wedge e^4) + \frac{(-a_5^2 + a_4^2)}{2 a_4} e^1 \wedge e^4 + a_4 e^2 \wedge e^3,\\
d e^2 = - \frac{a_5^2}{a_4} e^1 \wedge e^3 + a_5 (e^1 \wedge e^4 + e^2 \wedge e^3) +\frac{(a_5^2 -
a_4^2)}{2 a_4} e^2 \wedge e^4,\\
d e^3 = \frac{(a_4^2 + a_5^2)} {a_4} e^3 \wedge e^4, d e^4 = 0,
\end{array} \right.
$$
(with $\rho(e_1,e_1) = 0$ and $\rho(e_3,e_3) = - \frac{3}{2}  \frac{(a_4^2+a_5^2)^2}{a_4^2}$)
$$
17) \left\{
\begin{array} {l}
d e^1 = a_5 (e^1 \wedge e^3 + e^2 \wedge e^4) + \frac{ (a_5^2 - a_4^2)}{2 a_4} e^1 \wedge e^4 + a_4 e^2 \wedge e^3,\\
d e^2 = - \frac{a_5^2}{a_4} e^1 \wedge e^3 + a_5 (e^1 \wedge e^4 - e^2 \wedge e^3) +  \frac{(-a_5^2 +
a_4^2)}{2 a_4} e^2 \wedge e^4, \\
d e^3 = \frac{(-a_4^2 - a_5^2)} {a_4} e^3 \wedge e^4, d e^4 = 0,
\end{array} \right.
$$
(with $\rho(e_1,e_1) = 0$ and $\rho(e_3,e_3) = - \frac{3}{2}  \frac{(a_4^2+a_5^2)^2}{a_4^2}$).
\newline
If one puts $b_5 = - \frac{1}{2} b_2$ in the structure equations  5) one gets the structure
equations 11).  Moreover, the Lie group with structure equations 14) is isomorphic  to 5)
considering the new basis
$$
f^1 = -e^2, f^2= e^1, f^3 = e^3, f^4 = e^4.
$$
If one put $a_5 = 0$ in the structure equations $16)$ one gets $5)$. Then $5)$ is a particular case
of $16)$. Finally  $17)$ is isomorphic  to $16)$ considering the new basis 
$$
f^1 = e^1, f^2 = - e^2, f^3 = -e^3, f^4 = -e^4.
$$
and the respective metrics are isometric. 
Thus $G$ must have structure equations $16)$.
\end{proof}

\bigskip

\begin{remark} {\rm One has that $$\mbox {trace}
(ad_{e_4}) = \frac{t^2 + s^2} {t} \neq 0.$$ Therefore $G$ is not unimodular and it  does not admit any compact quotient
\cite{Mi}.}
\end{remark}

\medskip

\begin{remark} {\rm For any $t \neq 0$ and  $s \in \R$ the metric $g$ of Theorem 2.1 is never Einstein. Then, there exist no
left-invariant  Einstein, strictly almost K\"ahler metrics on real 4-dimensional Lie groups. \newline
If one investigates left-invariant (non flat)   K\"ahler-Einstein metrics on 4-dimensional Lie groups, by
using the proof of the above theorem and then by imposing the conditions
$\eqref{eq1},\eqref{eq2},\eqref{eq3},\eqref{eqintegrability}$ and
$$
\rho(e_1, e_1) = \rho(e_3, e_3),
$$
one gets the following cases for the relations between the non zero parameters: \newline
{\sl case 2)} with $b_5 \neq 0$: the
corresponding Lie group has structure equations:
$$
\begin{array} {l}
d e^1 =b_5 e^1 \wedge e^4 + a_5 e^2 \wedge e^4, d e^2 = - a_5 e^1 \wedge e^4
+ b_5 e^2 \wedge e^4,\\
d e^3 = 2 b_5 (e^1 \wedge e^2 + e^3 \wedge e^4), d e^4 = 0
\end{array}
$$
and  $\rho(e_1, e_1) = - 6 b_5^2$. \newline
{\sl case 3)} with $c_6 = \frac{a_4^3 + 3 a_6^2 a_4} {a_6^2 - a_4^2}$: the
corresponding Lie group has structure equations:
$$
\begin{array} {l}
d e^1 = \frac{a_6^3 + 3 a_4^2 a_6}  {a_6^2 - a_4^2}  e^1 \wedge e^2
-a_4 e^1 \wedge e^4  + a_4 e^2 \wedge e^3 + a_6 e^3 \wedge e^4, d e^2 = 0,\\
d e^3 =-a_4 e^1 \wedge e^2  + a_6 (e^1 \wedge e^4 - e^2 \wedge e^3)+\frac{a_4^3 + 3 a_6^2 a_4}{a_6^2 - a_4^2} e^3 \wedge e^4,
d e^4 = 0,
\end{array}
$$
and $\rho(e_1,e_1) = - 2 \frac{(a_4^2+a_6^2)^3}{(a_4^2 - a_6^2)^2}$. \newline
{\sl case 4)} with $a_1 =0$: the
corresponding Lie group has structure equations:
$$
\begin{array} {l}
d e^1 =  a_4 (- e^1 \wedge e^4 +
e^2 \wedge e^3), d e^2 = 0,\\
d e^3 = - a_4 (e^1 \wedge e^2 + e^3 \wedge e^4), d e^4 = 0,
\end{array}
$$
and $\rho(e_1, e_1) = -2 a_4^2$. This case  is equivalent  to the case 3) imposing  $a_6 = 0$.
\newline
{\sl case 8)}:  the
corresponding Lie group has structure equations:
$$
\begin{array} {l}
d
e^1 =- b_4 (e^1 \wedge e^3 + e^2 \wedge e^4) + \frac{b_4^2}{b_5} (e^1 \wedge e^4 - e^2 \wedge
e^3),\\
d e^2 = b_5 (e^1 \wedge e^3 + e^2 \wedge e^4) - b_4 ( e^1 \wedge e^4 - e^2 \wedge e^3),\\
d e^3 = \frac{b_4^2 + b_5^2} {b_5} ( e^1 \wedge e^2 + e^3 \wedge e^4), d e^4 = 0,
\end{array}
$$
and  $\rho(e_1, e_1) = - \frac{2(b_5^2+b_4^2)^2}{b_5^2}$. \newline
{\sl case 9)} with $a_1^2 + b_1^2 = c_6^2$: the
corresponding Lie group has structure equations:
$$
\begin{array} {l}
d e^1 = a_1 e^1 \wedge e^2, d e^2 = b_1 e^1 \wedge e^2,\\
d e^3 = c_6 e^3 \wedge e^4, d e^4 = 0
\end{array}
$$
and $\rho(e_1, e_1) = - a_1^2- b_1^2$. \newline
{\sl case 12)} with $c_6 = \frac{b_2^3 + 3 b_6^2 b_2} {b_2^2 - b_6^2}$: the
corresponding Lie group has structure equations:
$$
\begin{array} {l}
d e^1 = 0, d e^2 = \frac{(-b_6^3 -3 b_2^2 b_6)}{b_2^2 - b_6^2} e^1 \wedge e^2 + b_2 ( e^1 \wedge e^3 + e^2 \wedge e^4) +
b_6 e^3
\wedge e^4,\\ 
d e^3 = b_2 e^1 \wedge e^2 + b_6 (e^1 \wedge e^3 + e^2 \wedge e^4)  + \frac{b_2^3 + 3 b_6^2 b_2} {b_2^2 - b_6^2} e^3 \wedge
e^4, d e^4 = 0,
\end{array}
$$
and  $\rho(e_1, e_1) = - 2 \frac{(b_2^2+b_6^2)^3}{(b_2^2-b_6^2)^2}$. This case is equivalent to the case 8)
taking the new basis $\{ f^1 = - e^2, f^2 = e^1, f^3 = e^3, f^4 = e^4 \}$ and $b_2 = - a_4, b_6 = a_6$.
\newline
{\sl case 13)} with $c_6 = 0$:
the
corresponding Lie group has structure equations:
$$
\begin{array} {l}
d e^1 = a_6 ( e^1 \wedge e^2 + e^3 \wedge e^4), d e^2 = b_6 ( e^1 \wedge e^2 + e^3 \wedge e^4),\\
d e^3 = b_6 (e^1 \wedge e^3 + e^2 \wedge e^4) + a_6 (e^1 \wedge e^4 - e^2 \wedge e^3), d e^4 = 0,
\end{array}
$$
and $\rho(e_1,e_1) = -2 a_6^2-2 b_6^2$. \newline
Comparing all the above  possibilities  with the list of non-isomorphic Lie algebras given in \cite{Je}  in the
case 2) one gets the family of non isomorphic Lie algebras, for distinct values of $t$, with structure equations 
$$
\begin{array} {l}
d f^1 = 0,\\
d f^2 = f^1 \wedge f^2 - t f^1 \wedge f^3,\\
d f^3 = t f^1 \wedge f^2 -  f^1 \wedge f^3,\\
d f^4 = d f^2 = - 2 f^1 \wedge f^4 + 2 f^2 \wedge f^3, \quad 0 \leq t < \infty,
\end{array}
$$
and as Riemannian space each of these is the hermitian hyperbolic space.
In the remaining cases one gets the direct sum of a 2-dimensional Lie algebra with itself and as Riemannian space the direct
product of a 2-dimensional solvable group manifold, of curvature $-1$, with itself.}
\end{remark}

 A classification of 4-dimensional solvable Lie
algebras appears for example in  \cite{Ve} and a description of all simply-connected, 4-dimensional solvable real Lie groups
which have commutators of dimension 3 is also given in \cite{O}, where a classification of such Lie groups admitting
a left-invariant complex structure is given. We can prove that all the Lie algebras with structure equations 16) are isomorphic
to the Lie algebra of type ${\mathfrak g}_{4,9} (\alpha)$ with structure equations
$$
\begin{array} {l}
d f^1 = 0,\\
d f^2 = ( 1 - \alpha) f^1 \wedge f^2,\\
d f^3 = -  f^1 \wedge f^3,\\
d f^4 = -\alpha f^1 \wedge f^4 - f^2 \wedge f^3.
\end{array}
$$
and $\alpha = \frac{1}{2}$. The Lie algebra ${\mathfrak g}_{4,9} (2)$  is one of the Lie algebras with 
anti-self-dual non conformally flat inner product \cite{DS}.

\begin{teo} Let $G$ be a 4-dimensional Lie group admitting an invariant strictly almost K\"ahler structure $(g, J, \Omega)$
such that the Ricci tensor is $J$-invariant then its Lie algebra ${\mathfrak g}$ is isomorphic to the Lie algebra
${\mathfrak g}_{4,9} (\frac{1}{2})$ (in the notation of \cite{Ve}) with structure equations
\begin{equation} \label{vergne}
\begin{array} {l}
d f^1 = 0,\\
d f^2 = \frac{1}{2} f^1 \wedge f^2,\\
d f^3 = -  f^1 \wedge f^3,\\
d f^4 = -\frac{1}{2} f^1 \wedge f^4 - f^2 \wedge f^3
\end{array}
\end{equation}
and  $g$ is homothetic  to the metric
$$h =(f^1)^2 + (f^2)^2 +
(f^3)^2 + (f^4)^2. $$
\end{teo}

\begin{proof}
By the previous  Theorem,   ${\mathfrak g}$  has structure equations 16) and thus  it
 is an abelian extension of
the Lie algebra of the real 3-dimensional Heisenberg Lie algebra. Indeed, one has that $G$ is 3-step solvable with 
the commutator ${\mathfrak g}^1 = [{\mathfrak g}, {\mathfrak g}]$  the real 3-dimensional
Heisenberg Lie algebra, spanned by  $\{ Z, X, Y \}$ with $Z$ central and $[X,Y] = Z$, where
$$
Z = a_4 e_1 + a_5 e_2, \quad X = e_3, \quad Y =e_2.
$$
One can remark that $ad_{e_4}$ is a derivation of ${\mathfrak g}^1$ such that
$$
\begin{array} {l}
ad_{e_4} (Z) = \eta Z,\\
ad_{e_4} (X) = a Z + \mu X + \delta Y,\\
ad_{e_4} (Y) = b Z + k X + \nu Y,
\end{array}
$$
with
$$
\eta = \frac{(a_4^2 + a_5^2)}{2 a_4} = \frac{1} {2} \mu = - \nu, \quad  \delta = a =k = 0, \quad b =
\frac{a_5} {a_4}.
$$
Then if one consider on $G$ the new  basis $\{ A, Z, X, Y \}$ (where $A =\frac{1} {\mu} e_4$), one has the following
non zero Lie brackets
$$
[A, Z] =  \frac{1}{2} Z, \quad  [A, X] =  X, \quad  [A, Y] =  \frac{b} {\mu} Z  - \frac{1} {2} Y, \quad [X, Y] = Z.
$$
If one considers
the new basis  $$\left \{ f_1 = A, f_2 =  \frac{1}{\sqrt{a_4^2 + a_5^2}} (-Y +  \frac{b} {\mu} Z), f_3 = \frac{1}{\mu} X,
f_4 =
\frac{a_4} {(a_4^2 + a_5^2)^{3/2}}  Z
\right \}$$ (which is orthogonal with respect to
$g$)
 one has
$$
[f_1, f_4] = \frac{1}{2} f_4, \quad [f_1, f_2] = - \frac{1} {2} f_2, \quad [f_1, f_3] =  f_3, \quad [f_2, f_3] =
f_4
$$
and then one obtains the Lie algebra ${\mathfrak g}_{4,9} (\frac{1}{2})$ with structure equations
\eqref{vergne}. Moreover, $({\mathfrak g}, g)$ is isometric to ${\mathfrak g}_{4,9} (\frac{1}{2})$ endowed with the metric
$$
\frac{(a_4^2 + a_5^2)^2}{a_4^2} \left( (f^1)^2 + (f^2)^2 + (f^3)^2 + (f^4)^2 \right),
$$
since $\vert \vert f_i \vert \vert^2 = \frac{1} {\mu^2} = \frac{a_4^2} {(a_4^2 + a_5^2)^2}$, for any $i= 1, 2, 3, 4$.
\end{proof}

\begin{remark} {\rm It is easy to check that the Lie algebra ${\mathfrak g}_{4,9} (\frac{1}{2})$ is isomorphic to the Lie
algebra of $\R^2 \ltimes {\mbox Sol}_2$.}
\end{remark}

\section{Curvature}

If $(M, g)$ is an oriented 4-dimensional Riemannian manifold, the action of the Hodge $*$ operator 
on the bundle of the 2-forms $\Lambda^2 M$ induces the decomposition
$$
\Lambda^2 M = \Lambda^+ M \oplus \Lambda^- M,
$$
into the subbundles of self-dual and anti-self-dual 2-forms. If one considers the Riemannian
curvature tensor $R$ as a symmetric endomorphism of $\Lambda^2 M$ one has the following $SO
(4)$-splitting
$$
R = \frac{s}{12} Id \vert_{\Lambda^2 M} + \tilde  \rho_0 + W^+ + W^-,
$$
where $s$ is the scalar curvature, $\tilde \rho_0$ is the Kulkarni-Nomizu extension of the
traceless Ricci tensor $\rho_0$ to an endomorphism of $\Lambda^2 M$ (anti-commuting with $*$) and
$W^{\pm}$ are the self-dual and anti-self-dual parts of the Weyl tensor $W$.

If $(M, g, J)$ is an almost Hermitian 4-manifold and $\Omega$ the corresponding K\"ahler form, the
action of $J$ gives rise to the following orthogonal splitting 
\begin{equation} \label{selfdual}
\Lambda^+ M = \R \Omega \oplus \lbrack \lbrack \Lambda^{0,2} M \rbrack \rbrack,
\end{equation}
where $\lbrack \lbrack \Lambda^{0,2} M \rbrack \rbrack$ denotes the bundle of the bundle of the
$J$-anti-invariant real 2-forms. Then, the vector bundle ${\mathcal W}^+ = S^2_0 (\Lambda^+ M)$ of
the symmetric traceless endomorphisms of $\Lambda^+ M$ decomposes into the sum of the three
sub-bundles ${\mathcal W}^+_i$, $i = 1, 2, 3$ \cite{TV}. More precisely: \newline
${\mathcal W}_1^+$ is the trivial line bundle generated by the element $\frac{1}{8} \Omega \otimes
\Omega - \frac{1}{12} Id \vert_{\Lambda^+ M}$; \newline
${\mathcal W}_2^+ = \lbrack \lbrack \Lambda^{0,2} M \rbrack \rbrack$ is the subbundle of elements
which exchange the two factor in \eqref{selfdual}; \newline
${\mathcal W}_3^+ = S^2_0(\lbrack \lbrack \Lambda^{0,2} M \rbrack \rbrack)$ is the subbundle of
elements preserving the splitting \eqref{selfdual} and acting trivially on the first factor $\R
\Omega$.

Moreover
$$
\tilde {\rho_0} = \tilde {\rho_0}^{\mbox{inv}} + \tilde {\rho_0}^{\mbox{anti}},
$$
where ${\tilde \rho_0}^{\mbox{inv}}$ and ${\tilde \rho_0}^{\mbox{anti}}$ are the Kulkarni-Nomizu extensions of the
$J$-invariant and $J$-anti-invariant parts of the traceless Ricci tensor.

By \cite{AAD} any strictly almost K\"ahler 4-manifold whose curvature satisfies the three conditions
\newline 
(i) $\tilde{\rho_0}^{\mbox {anti}} = 0$, \hskip 20 pt (ii) $W_2^+ = 0$,  \hskip 20 pt
(iii) $W_3^+ = 0$ \newline
is locally isometric to the (unique)
4-dimensional proper 3-symmetric space. 

The relations (i),(ii),(iii) are closely related to the conditions defined by A. Gray \cite{Gr}
$$
\begin{array} {l}
(G1) \quad R_{XYZW} = R_{XYJZJW},\\
(G2) \quad R_{XYZW} = R_{JXJYZW} + R_{JXYJZW} + R_{JXYZJW},\\
(G3) \quad  R_{XYZW} = R_{JXJYJZJW}.
\end{array}
$$
Indeed, by \cite{AAD} an almost Hermitian 4-manifold $(M, g, J)$ satisfies the property $(G3)$ if and only
if the Ricci tensor is $J$-invariant and $W_2^+ = 0$. It satisfies $(G2)$ if, in addition, $W_3^+ = 0$. By \cite{AAD}, a
complete, simply connected strictly almost K\"ahler 4-manifold which satisfies the condition $(G2)$  is isometric to the
proper 3-symmetric space.

By Theorem 2.1 if $G$ is a simply connected real 4-dimensional Lie group with an invariant strictly almost K\" ahler
structure $(g, J, \Omega)$ such that the Ricci tensor is $J$-invariant, then the non zero
components of the Riemannian curvature
$R$  in terms of the orthonormal basis
$\{ e_1,
\ldots, e_4 \}$ are given by
$$
\begin{array} {l}
R_{1212} = R_{1234} = R_{3412} = \frac {1}{2} \frac{(t^2+ s^2)^2}{t^2},\\
R_{1313} =-R_{1324} = R_{1414} = R_{1423} = R_{2314} = R_{2323} = - R_{2413} = R_{2424} =- \frac {1}{4} \frac{(t^2+
s^2)^2}{t^2}, \\ R_{3434} =  - \frac{(t^2+ s^2)^2}{t^2}.
\end{array}
$$
By direct computation, one has that $W_2^+ = 0$, since 
$ \rho^* (R - L_3 R) = 0$, where
$$
\begin{array} {l}
\rho^* (R) (x,y) = \sum_{i = 1}^4 R(x,e_i,J y, J e_i),\\
(L_3 R) (x,y,z,w) = R(Jx, Jy, J z, J w)
\end{array}
$$
and $W_3^+ = 0$ since $\frac{1}{4} (I - L_2) (I + L_3) (R) = 0$, where
$$
(L_2 R) (x,y,z,w) = \frac{1}{2} [R(x,y,z,w) + R(Jx,Jy,z,w) + R(Jx,y,Jz,w) + R(Jx,y, z,Jw)].
$$
Then the Lie group $G$ satisfies the condition $(G2)$. As a consequence,

\begin{corol} If $G$ is a simply connected real 4-dimensional Lie group with an invariant strictly
almost K\"ahler structure $(g, J, \Omega)$ such that the Ricci tensor is $J$-invariant, then $G$ satisfies the condition
$(G2)$ and (up to homothety) is isometric to the (unique) 4-dimensional proper 3-symmetric space.
\end{corol}

\smallskip

{\bf Acknowledgements.} The author is grateful to Vestislav Apostolov, Sergio Garbiero and Simon Salamon for useful
conversations and remarks.

\bigbreak 
\renewcommand{\thebibliography}{\list{\arabic{enumi}.\hfil}
{\settowidth\labelwidth{18pt}\leftmargin\labelwidth\advance
\leftmargin\labelsep\usecounter{enumi}}\def\newblock{\hskip.05em}
\sloppy\sfcode`\.=1000\relax}\newcommand{\bi}{\vspace{-3pt}\bibitem}

\end{document}